
\input amstex

\documentstyle{amsppt}
\magnification=\magstep 1

\NoRunningHeads

\define\rank{\operatorname{rank}}
\define\Sym{\operatorname{Sym}}
\NoBlackBoxes

\topmatter

\title
Refined intersection products and limiting linear subspaces of
hypersurfaces
\endtitle{}

\author Xian Wu
\endauthor{}

\address
Department of Mathematics,
University of South Carolina,
Columbia, SC 29208
\endaddress{}

\email
wux{\@}milo.math.scarolina.edu
\endemail

\endtopmatter{}
\document

\smallpagebreak
\heading 0. Introduction
\endheading{}
\smallpagebreak

Let $X$ be a hypersurface of degree $d$ in $\Bbb P^n$
and $F_X$ be the scheme of $\Bbb P^r$'s contained in $X$.
If $X$ is generic,
then $F_X$ will have the expected dimension (or empty)
and its class is given by the top Chern class of
the vector bundle ${\Sym}^dU^*$,
where $U$ is the universal subbundle on
the Grassmannian $G(r+1,n+1)$.

Things become more interesting when $X$ degenerates.
For example,
when we deform a generic $X$ into a degenerate $X_0$,
the dimension of $F_X$ can jump.
In this case,
however,
there is a subscheme $F_{\text {\, lim}}$ of $F_{X_0}$
which consists of limiting
$\Bbb P^r$'s in $X_0$ with respect to the deformation.
For a generic deformation,
$F_{\text {\, lim}}$ will have the expected dimension
and its class must also be determined by the top
Chern class of ${\Sym}^dU^*$.
This subscheme $F_{\text {\, lim}}$
of limiting $\Bbb P^r$'s will be the main
object of study in this paper.

Geometrically,
there are two basic questions concerning
$F_{\text {\, lim}}$ that arise naturally.
Let us take a simple example.
We all know that there are $27$ lines in a generic cubic surface.
If we degenerate the surface into the union of
a plane and a quadric,
then there are infinitely many lines in the union.
Which $27$ lines are the limiting ones
and how many of them are in the plane
and how many of them are in the quadric?
In general,
we would like to ask the following questions:
How can we identify limiting $\Bbb P^r$'s
to characterize $F_{\text {\, lim}}$?
If $F_{X_0}$ has more than one component,
what is the local distribution of limiting
$\Bbb P^r$'s to those components?
Such problems in the case of $r=1$ are investigated in [W].
The main tools used there are infinitesimal deformation theory
and the theory of Chern classes.
In particular,
the answer to the example of the cubic surface above is very simple.
The intersection of the plane, the quadric,
and a cubic (which depends on the deformation)
consists of six points.
A line in the plane is a limiting line if and only if it passes through
any two of those six points.
That gives us $15$ limiting lines in the plane.
On the other hand,
a line in the quadric is a limit if and only if it passes through
any one of those points.
Since for every point in a generic quadric
there are two lines that pass through it,
we see that there are $12$ limiting lines in the quadric.

One of the main goals of this paper is to answer the questions above
for arbitrary $r$.
In particular,
we give a simple but concrete geometric
description of $F_{\text {\, lim}}$
and explicit formulas to compute the local distribution
for some generic degenerations.
Please see Summary/Theorem 4.1 in Section 4
for the precise statement.

Our formulas (Theorem 3.1 and Proposition 5.3)
for the local distribution are computable by
using standard techniques in the theory of Chern classes
and Schubert calculus.
We have included some examples in Section 5.
Let us just mention one (Example 4) here since it seems to
be a bit surprising.
We see that there are $321,489$
$\Bbb P^3$'s on a generic
cubic hypersurface in $\Bbb P^8$.
When the cubic degenerates into the union of a quadric and a plane,
computations from our formulas reveal that
all limiting $\Bbb P^3$'s
are in the plane and none of them is in the quadric!
In particular,
this implies that $U^*$ is not a positive bundle.

Since the total class of $F_{\text {\, lim}}$
is always equal to the top Chern class of
${\Sym}^dU^*$,
as a consequence of our main theorems,
we get a family of identities
involving the Chern classes and the Segre classes
of various symmetric powers of $U^*$
(Corollary 5.1 and Corollary 5.7).
We think that those formulas should hold for any vector bundle
of rank $r+1$ over general manifolds.
In fact,
they should be identities of symmetric polynomials in $r+1$ variables,
as we have seen in the case of $r=1$ in [W].
It will be interesting to see if there are other possible applications
of such formulas.
Those identities can also be found in Section 5.

The case $r>1$ brings some fundamental changes from $r=1$.
The first,
there seem to be some technical difficulties to apply
the infinitesimal deformation theory in the way used in [W].
However, the main difference is that
we can no longer ignore the scheme of $\Bbb P^r$'s contained
in the singular locus of $X$ if $r>1$.
This is a common difficulty that one often faces in
doing computations
in intersection theory.
As a results,
we see that the formulas obtained for general $r$
is much more complicated than those for lines
and they are not generalizations of the formulas in [W].
On the other hand,
we believe that such degeneration problems
are not isolated ones
and there should be an intrinsic treatment in general.
As it turns out,
the theory of refined intersection products of Fulton-MacPherson
provides a natural way to study our geometric problems.
The idea is very simple.
Any hypersurface $X$ induces a section $s_X$ of ${\Sym}^dU^*$.
The zero-scheme of $s_X$ is nothing but $F_X$.
Assuming that $F_X$ is not empty,
we can then consider the fiber square
$$
\CD
F_X @>>> G(r+1,n+1)\\
@VVV @VVs_XV\\
G(r+1,n+1) @>i>> {\Sym}^dU^*
\endCD{}
$$
where $i$ is the zero-section imbedding.
{}From this we can construct a class,
called ${R}_{X}$ in this paper,
by the refined intersection product of Fulton-MacPherson
induced from the above fiber square
[F, Chapter 6] [FM].
By well-known facts about refined intersection products and
dynamic intersections,
this class has the following interesting properties (Theorem 1.3).

(1) It is defined pointwise from data on $X$
and hence is independent of deformations.

(2) Despite (1),
its image in the Chow ring of $G(r+1,n+1)$
is always equal to the top Chern class of ${\Sym}^dU^*$
and is hence independent of the choice of $X$.

(3) Geometrically,
it equals to the class of $F_X$
if $X$ is generic
and otherwise is equal to the class of the subscheme $F_{\text {\, lim}}$
of limiting $\Bbb P^r$'s
with respect to a generic deformation.

(4) Finding the local distribution of limiting $\Bbb P^r$'s
is equivalent to finding the canonical decomposition of
${R}_X$ in terms of its distinguished varieties.

What makes this formulation more useful is the fact that
the above construction and properties are valid for much more general
settings under very modest hypotheses.
However,
unless their distinguished varieties are connected components,
which is only possible for $r=1$ in the present case,
such intersection products
and their canonical decompositions are in general difficult to compute.
So the construction above does not give us explicit answers directly.
Nevertheless,
such a formulation sets up the stage for
a nice interplay between our
geometric problems and intersection theory.
In particular,
it allows us to apply a powerful theorem of Fulton and Lazarsfeld [F] [La]
on infinitesimal intersection class
to conclude that $F_{\text {\, lim}}$ is determined by
the infinitesimal data (Theorem 2.2 and Corollary 2.3).
This is one of the key steps in getting our main geometric results.
On the other hand,
by applying the geometric results back to intersection theory,
we are able to find the equivalences of the distinguished varieties and
hence the canonical decomposition for refined intersection
product ${R}_{X}$ constructed above (Summary/Theorem 4.1').

While the problem of understanding limiting subschemes
arises naturally from studies of
algebraic cycles using degeneration techniques,
the formulation in the refined intersection product
suggests that it may have interesting connections with other
problems of a different nature.
We have been seeing recently with great excitement the
dramatic interplay between string theory and algebraic geometry.
In particular,
the development has led to some amazing computations
[COGP] [M1] [M2] which
give striking
predictions for the numbers of rational curves in certain Calabi-Yau
threefolds and,
more recently [GMP],
in higher dimensional cases.
The Calabi-Yau manifolds in questions are often special types which
contain families of rational curves.
So one could explain that what those authors have really computed are
the numbers of limiting rational curves.
Interestingly,
no deformations of those Calabi-Yau manifolds
are presented in the method
(there are only deformations of their ``mirrors'').
On the other hand,
the definition of ${R}_X$
is also purely pointwise and it can be easily generalized
to study the limiting rational curves.
It is our hope that
the basic idea of this paper may lead to more useful results.

This paper is organized as follows.
We start our study in Section 1
by giving the construction of ${R}_X$
and interpreting our geometric problems
in terms of intersection theory (Theorem 1.3).
This is mainly a direct application of some well-known facts
about dynamic intersection theory and
the refined intersection products of Fulton-MacPherson.

In Section 2,
we study the subscheme $F_{\text {inf}}$ of infinitesimal
limiting $\Bbb P^r$'s.
A main tool here is the infinitesimal intersection class
as defined in [F, Chapter 11].
By a theorem of Fulton and Lazarsfeld [F, Chapter 11] [La],
this infinitesimal intersection class refines the class of
$F_{\text {\, lim}}$.
Moreover,
for a generic deformation,
$F_{\text {inf}}$ will have the expected dimension and its
class is exactly
equal to the above infinitesimal intersection class.
In other words,
limiting $\Bbb P^r$'s are determined by
the infinitesimal data in that case (Corollary 2.3).
With these general facts in place,
we proceed with some explicit computations.
For this,
we consider a generic deformation $\Cal D$
which degenerates $X$ into the union of two hypersurfaces.
Up to first order,
we may write such a deformation as
$$
\Cal D=\{X_s=KL+sD \}_{s\in \Delta},
$$
where $K$, $L$, and $D$
are generic hypersurfaces of degree $k$, $l$, and $d=k+l$,
respectively.
We will use the same letter to denote a hypersurface and its
equation.
Let
$$
\sigma_K=\{\Bbb P^r \in G(r+1,n+1)|
\  \Bbb P^r \subset K, \Bbb P^r \cap L \subset D\}
$$
and
$$
\sigma_L=\{\Bbb P^r \in G(r+1,n+1)|
\  \Bbb P^r \subset L, \Bbb P^r \cap K \subset D\}.
$$
The fact (Proposition 2.5) is that
$\sigma_K$
and $\sigma_L$
represent the subschemes of infinitesimal limiting $\Bbb P^r$'s
and hence limiting $\Bbb P^r$'s
in $K$ and $L$,
respectively.
This is done by using the method of deforming ideals adopted from
Katz's earlier work [K2].

Turning to Section 3,
we will compute the class of $\sigma_K$ and the class of $\sigma_L$.
This has been done in [W] for the case of $r=1$,
since the geometric conditions
which define $\sigma_K$ and $\sigma_L$ are the same
for $r=1$ and $r>1$.
However,
the method used in [W] cannot be generalized here.
The main difficulty is that,
while $\Bbb P^1$'s contained in both $K$ and $L$
can basically be ignored,
the same is no longer true for $\Bbb P^r$'s if $r$ is greater than $1$.
In terms of intersection theory,
that means while in the case of $r=1$
we can treat the distinguished varieties
as though they were connected components,
we are now dealing with the excess intersections in a
much more complicated way.
In fact,
a straightforward generalization of the formulas
in [W] gives obviously wrong answers.
Therefore,
some correction terms must be added.
Fortunately,
nowadays we have the power of the modern technology.
After many ``experiments'' on computers,
we were able to guess correct formulas (Theorem 3.1).
The actual proof of those formulas
follows naturally once we know what to prove.

Section 4 is used primarily to give
a summary of the main results obtained in the first
three sections for the convenience of readers.

In Section 5,
we give some examples and applications
of our main results.
In particular,
as mentioned before,
we obtain a family of identities
in the characteristic classes of
various bundles involving $U^*$
(Corollary 5.1 and Corollary 5.7).
Those identities seems to be complicated enough to be non-trivial and,
on the other hand,
to be in a form nice enough to make one suspect that
they may have other meanings and applications.
Most examples are computed by
using a Maple package called ``schubert'' [KS]
written by Katz and Str{\o}mme.

Finally,
for the convenience of interested readers,
we have included schubert code for calculations of two examples
in the appendix.
We want to thank Katz and Str{\o}mme
for making this wonderful Maple package available
as well as for detailed help on actually using it.

The author is indebted to S. Katz for many valuable discussions.
The author also wants to express his thanks to W. Fulton
for his encouragement and advice.
A part of the work was done while the author was a participant in
the NSF Regional Geometry Institute hold at Amherst College
during the summer, 1992.
The author wants to thank the organizers and the NSF for their generous
support.

\smallpagebreak
\heading 1. Refined intersection products and
limiting linear subspaces
\endheading{}
\smallpagebreak

In this section we will lay out the ground work by applying
some well-known facts from intersection theory to our geometric problems.
Although the statements here
are given in forms appropriate for our present setup,
most of them can be easily generalized to other settings.
The main reference for this section is Fulton's book [F].

Let $G(r+1,n+1)$ be the Grassmannian of
$\Bbb P^r$ in $\Bbb P^n$
and $\Bbb P^N$ be the Hilbert scheme
of hypersurfaces of degree $d$ in $\Bbb P^n$,
where
$$
N=\dim(\Bbb P(H^0(\Cal O_{\Bbb P^n}(d))))={n+d \choose d}-1.
$$
Throughout this paper,
we will use the same notation for a point of $G(r+1,n+1)$ and
its corresponding $\Bbb P^r$ in $\Bbb P^n$.
Similarly,
given an $X$ in $\Bbb P^N$,
we will use same $X$ to denote the corresponding hypersurface
in $\Bbb P^n$ and its homogeneous equation.
Consider the correspondence
$$
\CD
\Cal Y @>{p_1}>> G(r+1,n+1) \\
@VVp_2V  @.\\
\Bbb P^N
\endCD{}
\tag 1.1
$$
where
$$
\Cal Y = \{(\Bbb P^r,X)\ |\  \Bbb P^r \subset X\}
\ \subset \ G(r+1,n+1)\times\Bbb P^N.
$$
Given an $X$ in $\Bbb P^N$,
let $F_X$ be the scheme of $\Bbb P^r$'s in $X$,
that is,
$$
F_X=p_1(p_2^{-1}(X)).
$$
The fiber over any $\Bbb P^r$ in $G(r+1,n+1)$ is
$\Bbb P(H^0(\Cal I_{\Bbb P^r/\Bbb P^n}(d)))$.
{}From the standard exact sequence
$$
0 \to \Cal I_{\Bbb P^r/\Bbb P^n}(d)
\to \Cal O_{\Bbb P^n}(d)
\to \Cal O_{\Bbb P^r}(d) \to 0,
$$
it is easy to see that $\Cal Y$ is irreducible and smooth of
dimension
$$
(r+1)(n-r)+N-{r+d \choose d}.
$$
Hence,
for a generic $X$ in $\Bbb P^N$,
$F_X$  is either empty or smooth of dimension $m$,
where
$$
m=(r+1)(n-r)-{r+d \choose d}.\tag 1.2
$$
We will assume for a moment that $F_X$ is not empty and,
in particular,
$m$ is non-negative.
For any special $X_0$ in $\Bbb P^N$ such that the dimension of $F_{X_0}$
is greater than $m$,
we will consider limiting $\Bbb P^r$'s in $X_0$.
The following consideration will be local.
For any one-parameter deformation $\Cal D$ of $X_0$ in $\Bbb P^N$,
set
$$
\Cal Y_{\Cal D}=\overline{p_2^{-1}(\Cal D-X_0)}.
$$
We define the scheme of limiting $\Bbb P^r$'s in $X_0$ with respect
to $\Cal D$ to be
$$
F_{\text {\, lim}}^{\Cal D}
=\lim_{ X {\overset\Cal D\to \to} X_0} F_X
=p_1(\Cal Y_{\Cal D}\cap p_2^{-1}(X_0)).
$$
{}From the definition,
$F_{\text {\, lim}}^{\Cal D}$ is a subscheme of $F_{X_0}$.
If $\Cal D$ is generic,
then $p_2^{-1}(X)$ has dimension $m$
for $X$ in $\Cal D-X_0$.
Therefore,
$F_{\text {\, lim}}^{\Cal D}$
also has dimension $m$.
In this case,
we will drop the superscript and denote it by $F_{\text {\, lim}}$.
Let $U$ be the universal subbundle on $G(r+1,n+1)$
and ${\Sym}^dU^*$ be
the $d^{th}$ symmetric power of the dual of $U$.
Furthermore,
let $c({\Sym}^dU^*)$ be the total Chern class of ${\Sym}^dU^*$
and $s(F_X, G)$ be the Segre class of $F_X$ in $G(r+1,n+1)$.

\proclaim{Theorem 1.3}
For any $X$ in $\Bbb P^N$,
there is a well-defined class ${R}_{X}$
in the Chow ring $A_m(F_X)$ of $F_X$
given by
$$
\left\{
\aligned
{R}_{X}=&\{c({\Sym}^dU^*)\cap s(F_X,G)\}_m,\\
{R}_{X}=&0, \endaligned{}
\quad \aligned \text{if} \ F_X &\not = \emptyset,\\
\text {if} \ F_X&=\emptyset.
\endaligned{}
\right.
\tag{1.4}
$$
such that:

(1) its image in $A_m(G(r+1,n+1))$ is independent of $X$ and in fact
$$
j_*({R}_{X})=c_{top}({\Sym}^dU^*)\cap G(r+1,n+1),
$$
where $j$ is the inclusion of $F_X$ in $G(r+1,n+1)$.

(2) For a generic $X$ in $\Bbb P^N$,
${R}_{X}$ is equal to the class of $F_X$,
that is,
$$
{R}_{X}=[F_{X}].
$$

(3) For a special $X=X_0$ such that
$F_{X_0}$ has a dimension greater than $m$,
${R}_{X_0}$ represents the subscheme of limiting $P^r$'s in $X_0$
with respect to a generic deformation $\Cal D$.
In other words,
$$
{R}_{X_0}=[F_{\text {\, lim}}].
$$

(4) If $F_i$ is a component of $F_{X_0}$,
then the part of ${R}_{X_0}$ supported on
$F_i$ represents the part of limiting
$\Bbb P^r$'s contained in $F_i$.
Equivalently,
$$
{R}_{X_0}^{F_i}=[F_{\text {\, lim}}\cap F_i].
$$
\endproclaim{}

\demo{Proof of Theorem 1.3}
Theorem 1.3 is an immediate consequence of well-known facts from
intersection theory once a connection between them is established.
We will use constructions and facts from
Fulton's book [F],
mainly on refined intersection products and dynamic intersections.
Please refer to [F]
for further details.

For any $X$ in $P^N$,
$X$ induces a section $s_X$ of ${\Sym}^dU^*$.
In fact,
${\Sym}^dU^*$ is generated by such sections.
This follows from the fact that the natural map
$$
H^0(\Cal O_{\Bbb P^n}(d))_{G(r+1,n+1)} \to {\Sym}^dU^*
$$
is surjective.
The zero-scheme $Z(s_X)$ of $s_X$ consists of those $\Bbb P^r$'s
which are contained in $X$.
In other words,
$Z(s_X)$ is equal to $F_X$.
If $F_X$ is empty,
then using the splitting principle it is easy to see that
the top Chern class of ${\Sym}^dU^*$ is zero.
We are therefore done in this case.
On the other hand,
if $F_X$ is not empty,
we can consider the fiber square
$$
\CD
F_X @>>> G(r+1,n+1)\\
@VVV @VVs_XV\\
G(r+1,n+1) @>i>> {\Sym}^dU^*
\endCD{}
\tag 1.5
$$
where $i$ is the zero-section imbedding.
We define ${R}_{X}$ to be the refined intersection product induced from
fiber square (1.5) as given in [F, Chapter 6].
In other words,
$$
{R}_{X}=i^!([G(r+1,n+1)])=i^*(C_{F_X}),
$$
where $C_{F_X}$ is the normal cone to $F_X$ in $G(r+1,n+1)$,
$i^!$ is the refined Gysin homomorphism induced from (1.5),
and $i^*$ is the Gysin homomorphism induced from $i$.
Since the normal bundle to
$G(r+1,n+1)$ in ${\Sym}^dU^*$ is ${\Sym}^dU^*$ itself,
${R}_{X}$ therefore is a well-defined class in $A_m(F_X)$
given by
$$
{R}_{X}=\{c({\Sym}^dU^*)\cap s(F_X,G)\}_m
$$
and its image in $A_*(G(r+1,n+1))$ is equal to
$$
c_{top}({\Sym}^dU^*)\cap G(r+1,n+1).
$$
To see claim (2),
recall that ${\Sym}^dU^*$ is generated by sections induced from
hypersurfaces of degree $d$ in $\Bbb P^n$.
Therefore,
for a generic $X$ in $\Bbb P^N$,
the induced section $s_X$ is transverse to
the zero-section.
Hence,
by the definition,
${R}_{X}$ is equal to the class of $F_X$.
(This also gives another way to see that $F_X$ is smooth
of dimension $m$.)
Claim (3) and Claim (4)
will now follow from the definition of $F_{\text {\, lim}}$
and general facts about dynamic intersections.
Taking a deformation $\Cal D$ and deforming fiber square (1.5)
accordingly,
we get the fiber square
$$
\CD
\Cal F @>>> G(r+1,n+1)\times \Cal D\\
@VVV @VVsV\\
G(r+1,n+1)\times \Cal D @>i\times id>> {\Sym}^dU^*\times \Cal D
\endCD{}
\tag 1.6
$$
where
$$
s(*,X)=(s_X(*),X).
$$
Define ${\Cal {R}}$ to be the refined intersection product
in $A_{m+1}(\Cal F)$
induced from fiber square (1.6).
Notice that $\Cal F$ is a scheme over $D$ with the fiber over $X$
being $F_X$.
Since the rational equivalence is preserved when passing
from a family to its fibers,
${\Cal {R}}$
defines a family of classes over $\Cal D$ with the fiber over
$X$ being ${\Cal {R}}_X$.
Next,
let us recall the definition of the limit of ${\Cal {R}}_X$.
Assume
$$
{\Cal {R}} =\sum_{i}n_i[\Cal V_i],
$$
where $\Cal V_i$ are subvarieties of $\Cal F$.
Let
$$
\Cal V_i^*=\Cal V_i - F_{X_0}.
$$
We define the limit of ${\Cal {R}}_X$ to be
$$
\lim_{ X {\overset\Cal D\to \to} X_0} {\Cal {R}}_X
=\sum_i n_i[\overline{\Cal V_i^*}\cap F_{X_0}].
$$
It is easy to check that
this limit is a well-defined class.
An important fact in dynamic intersections [F, Chapter 11]
is that
$$
\lim_{ X {\overset\Cal D\to \to} X_0} {\Cal {R}}_X ={\Cal {R}}_{X_0}.
$$
On the other hand,
fiber square (1.6) can be considered as a family of fiber squares
over $\Cal D$ with fibers over $X$ being fiber square (1.5).
{}From this,
we see that ${\Cal {R}}_X$ is nothing but ${R}_{X}$.
Therefore,
for a generic deformation $\Cal D$,
$$
\aligned
{R}_{X_{0}}=&{\Cal {R}}_{X_0}\\
=&\lim_{ X {\overset\Cal D\to \to} X_0} {\Cal {R}}_X\\
=&\lim_{ X {\overset\Cal D\to \to} X_0} {R}_{X}\\
=&\lim_{ X {\overset\Cal D\to \to} X_0} [F_X]\\
=&[\lim_{ X {\overset\Cal D\to \to} X_0} F_X]\\
=&[F_{\text {\, lim}}].
\endaligned{}
$$
{}From the construction,
it is clear that Claim (4) also follows.
\qed
\enddemo{}

Since $F_X$ gives a positive class if it is not empty,
we have the following corollary.

\proclaim{Corollary 1.7}
For a generic $X$,
the scheme $F_X$ is empty if and only if $c_{top}({\Sym}^dU^*)$
is zero.
\endproclaim{}

\example{Remark}
The scheme $F_X$ can be empty even if $m$
is non-negative.
For example,
it is well-known [GH, Chapter 6]
that a smooth quadric in $\Bbb P^{n}$
contains no $\Bbb P^r$ if $r$ is strictly greater than
$(n-1)/2$.
On the other hand,
for $d=2$ in (1.2),
$m$ is non-negative if $r$ is less than or equal to $2(n-1)/3$.
\endexample{}

{}From Theorem 1.3,
it is clear now that the next step is to understand ${R}_{X}$.
Of course,
the image of ${R}_{X}$ in $A_m(r+1,n+1)$
is always equal to the top Chern class of ${\Sym}^dU^*$.
Our geometric goal is to identify $F_{\text {\, lim}}$ in $F_{X}$
if the dimension of $F_X$ is greater than $m$.
In terms of ${R}_{X}$,
that means to compute its canonical decomposition.
For example,
$F_X$ will have $l$ components $F_{X_i}$,
$1\leq i \leq l$,
if we degenerate a generic hypersurface into a normal crossing
of $l$ generic hypersurfaces
$$
X=\bigcup_{i=1}^l X_i, \quad \sum_{i=1}^l\deg X_i=d.
$$
Those $F_{X_i}$ are all distinguished varieties of ${R}_{X}$
if we assume that they are all irreducible.
(This may not be the case in general,
but that is not important for our discussion.)
To find limiting $\Bbb P^r$'s in each $X_i$ is equivalent to
find the equivalence, or the contribution,
of corresponding distinguished variety to ${R}_{X}$.
Unfortunately,
unless those varieties are also connected components,
such a decomposition is in general difficult to compute.
In fact,
one cannot expect an explicit formula in such a general setting.
The following is how we will proceed.
First, we will see that the problem can be reduced to
a problem of the infinitesimal nature.
To make more explicit computations,
we then turn our attention to degenerations such that $X_0$
has two components and
identify infinitesimal limiting $\Bbb P^r$'s
in each component geometrically.
Finally,
the class of limiting $\Bbb P^r$'s in each component will be computed
from those geometric conditions.
Those steps will be carried out
in next two sections.

\smallpagebreak
\heading 2. Infinitesimal intersection classes and
infinitesimal limiting linear subspaces
\endheading{}
\smallpagebreak

In this section we will study the subscheme of
infinitesimal limiting $\Bbb P^r$'s.
Our considerations therefore will be all local.
We will use a theorem of Fulton and Lazarsfeld on infinitesimal
intersection classes.
This makes it possible to study the subscheme of limiting $\Bbb P^r$'s
from infinitesimal data.
The theorem of Fulton and Lazarsfeld implies that
this principle is true in much more general settings.
Please refer to [F, Chapter 11] for further details.

Let us first recall a construction of [F, Chapter 11].
For a given $X_0$ in $\Bbb P^N$ and a one-parameter deformation
$\Cal D$ of $X_0$,
we consider the imbedding $s$ as given in (1.6)
$$
G(r+1,n+1)\times \Cal D \ \overset s \to \longrightarrow
\ {\Sym}^dU^*\times \Cal D.
\tag 2.1
$$
Let $T$ be the tangent space to $\Cal D$ at $X_0$
and $T_G$ be the trivial line bundle $T\times G(r+1,n+1)$ on $G(r+1,n+1)$.
Let $N_{X_0}$ be the normal bundle of
$$
s_{X_0}(G(r+1,n+1)) \simeq s(G(r+1,n+1),X_0)
$$
in ${\Sym}^dU^*$.
Since $F_{X_0}$ is the zero-scheme of $s_{X_0}$,
we have that
$$
N_{X_0}|_{F_{X_0}}={\Sym}^dU^*|_{F_{X_0}}
$$
and we will call this restriction $N$.
{}From (2.1),
a basis $\partial / \partial t$ of $T$
will induce a section of $T_G$ and hence a section of $N_{X_0}$ by
the Kodaira-Spencer homomorphism.
Let $s_T$ be the restriction of the above section of $N_{X_0}$ on $N$.
Hence $s_T$ is a section of ${\Sym}^dU^*$ over $F_{X_0}$.
Let $C$ be the normal cone to $F_{X_0}$ in $G(r+1,n+1)$.
Notice that $C$ is a purely $(n-r)(r+1)$-dimensional subscheme of $N$.
Let
$$
F_{\text {inf}}=s_T^{-1}(C)
$$
be the subscheme of infinitesimal limiting $\Bbb P^r$'s in $F_{X_0}$.
The infinitesimal intersection class
${R}_{\text {inf}}$ is then defined to be
$$
{R}_{\text {inf}}=s_T^!(C) \ \in \ A_m(F_{\text {inf}}).
$$

We can now state a theorem of Fulton and Lazarsfeld [F] [La].

\proclaim{Theorem 2.2 (Fulton and Lazarsfeld)}
We have the following inclusions
$$
F_{\text {\, lim}} \hookrightarrow
F_{\text {inf}} \overset j \to \hookrightarrow
F_{X_0}.
$$
Moreover,
${R}_{\text {inf}}$ refines ${R}_{X_0}$,
that is,
$$
j_*({R}_{\text {inf}})={R}_{X_0}.
$$
\endproclaim{}

\demo {Proof}
Since $s_X$ is a regular imbedding of $G(r+1,n+1)$ in ${\Sym}^dU^*$,
the refined intersection product ${\Cal {R}}$ defined in the proof
of Theorem 1.3 can also be defined from the transform of
fiber square (1.6)
$$
\CD
\Cal F @>>> G(r+1,n+1)\times \Cal D\\
@VVV @VVi\times idV\\
G(r+1,n+1)\times \Cal D @>s>> {\Sym}^dU^*\times \Cal D
\endCD{}
$$
Similarly,
consider ${R}_{X}$ as the refined intersection product
induced from the transform of
fiber square (1.5).
Now apply Theorem 11.2 of [F].
\qed
\enddemo{}

If we take $\Cal D$ to be general enough,
$F_{\text {inf}}$ will have the expected dimension and
${R}_{\text {inf}}$ then represents the class of $F_{\text {inf}}$.

\proclaim{Corollary 2.3}
For a generic deformation $\Cal D$,
$$
[F_{\text {\, lim}}] = [F_{\text {inf}}] = {R}_{X_0}
$$
and the subscheme $F_{\text {\, lim}}$ of limiting
$\Bbb P^r$'s is hence determined by
the infinitesimal data.
\endproclaim{}

We will now make some explicit computations about $F_{\text {inf}}$.
For this,
we will restrict our attention to a simple type
of degenerations.
To be more precise,
we will consider
$$
\Cal D=\{X_s\}_{s\in \Delta}
$$
a semistable degeneration of generic hypersurfaces of degree $d$
in $\Bbb P^n$ such
that $X_0$ has two irreducible components.
As indicated in the introduction,
we will use the same letter to denote a hypersurface
and its equation.
Up to first order,
we may write that
$$
\Cal D =\{X_s=KL+sD\}_{s \in \Delta},\tag 2.4
$$
where $K$, $L$, and $D$ are generic hypersurfaces of degree $k$, $l$,
and $d=k+l$,
respectively.
We define
$$
\sigma_K = \{\Bbb P^r \in G(r+1,n+1)|
\  \Bbb P^r \subset K, \ \Bbb P^r \cap L \subset D\}
$$
and
$$
\sigma_L = \{\Bbb P^r \in G(r+1,n+1)|
\  \Bbb P^r \subset L, \ \Bbb P^r \cap K \subset D\}.
$$
Notice that $\sigma_K \cup \sigma_L$ is a subset of $F_{X_0}$
subject to the given geometric conditions.
The extra conditions are imposed exactly by the first order obstructions
to the deformation.
To be more precise, we have the following proposition.

\proclaim{Proposition 2.5}
For a generic degeneration (2.4),
$$
F_{\text {inf}} = \sigma_K \cup \sigma_L .
$$
\endproclaim{}

\demo{Proof of Proposition 2.5}
We will use the method of deforming ideals
to study the first order obstruction on $\Bbb P^r$'s
in $F_{X_0}$ as
Katz [K2] did in the case of quintic threefolds.
For any $\Bbb P^r$ in $F_{X_0}$,
we have that either
$\Bbb P^r$ is contained in $ K$
or
$\Bbb P^r$ is contained in $L$.
Assume that
$\Bbb P^r$ is contained in $K$ and it deforms with $X_s$.
Take a local trivialization of $\Cal D$
and a family
$$
\Bbb P_{s}^r \subset X_s,\quad \Bbb P_0^r= \Bbb P^r.
$$
Let $I_s$ be the ideal of $\Bbb P_{s}^r$ in $\Bbb P^n$
and
$$
I_s=\{f_1(s), f_2(s), \cdots, f_t(s)\}.
$$
Since $X_s$ is in $I_s$,
we have that
$$
X_s=\sum_{i=1}^t a_i(s)f_i(s)=KL+sD+{s}^2E+\cdots .
$$
Therefore,
we have infinitesimally that
$$
\left\{
\eqalign{KL&=\sum_{i=1}^t a_i(0)f_i(0),\cr
D &= \sum_{i=1}^t a_i(0)f_i'(0)+\sum_{i=1}^t a_i'(0)f_i(0).\cr}
\right .
\tag 2.6
$$
Since $\Bbb P^r$ is in $K$,
$$
K=\sum_{i=1}^t b_if_i(0).
$$
Therefore,
from the first equation in (2.6),
we may choose $K$ and $L$ such that
$$
a_i(0)=b_iL.
$$
Substituting this into the second equation in (2.6),
we get
$$
D=L\sum_{i=1}^t b_if_i'(0)+\sum_{i=1}^t a_i'(0)f_i(0).
$$
Since
$$
f_i(0)|_{\Bbb P^r}=0,
$$
we then have that
$$
D|_{\Bbb P^r}=L|_{\Bbb P^r}\left. \left(\sum_{i=1}^t b_if_i'(0)\right)
\right|_{\Bbb P^r}.
$$
The equation above implies that
$$
D|_{\Bbb P^r}(p)=0,\quad
\text{whatever} \ \ L|_{\Bbb P^r}(p)=0.
$$
This is just the definition of
$\sigma_K$.
\qed
\enddemo{}

Our next step will be
to compute the class of $\sigma_K$ and the class of $\sigma_L$
from the given geometric conditions.
For this,
we will turn to Section 3.

\smallpagebreak
\heading 3. Formulas for $\sigma_K$ and $\sigma_L$
\endheading{}
\smallpagebreak

Our goal in this section is to derive formulas for
the classes of $\sigma_K$ and $\sigma_L$
in the Chow ring $A_*(G(r+1,n+1))$ from their definitions.
That together with the results in Section 1 and Section 2
will enable us to compute the local distribution of limiting $\Bbb P^r$'s
and hence the canonical decomposition
of ${R}_{X_0}$.
To simplify the notation,
we will identify an element in $A^*(G(r+1,n+1))$
with its dual in $A_*(G(r+1,n+1))$.

Retaining the notation in previous sections,
we consider first the projective bundle
$$
\Bbb P({\Sym}^lU^*) @>\pi_l>> G(r+1,n+1).
$$
Let
$$
T_l=\Cal O(-1)
$$
be the tautological subbundle on $\Bbb P({\Sym}^lU^*)$.
Similarly,
let $T_k$ be the tautological subbundle on the projective bundle
$$
\Bbb P({\Sym}^kU^*) @>\pi_k>> G(r+1,n+1).
$$
For a given vector bundle $E$,
we will use $c_{top}(E)$ to denote the top Chern class of $E$.
Recall also that we have $d=k+l$.

\proclaim{Theorem 3.1}
Let $[\sigma_K]$ be the class of $\sigma_K$ in $A_*(G(r+1,n+1))$.
We have that
$$
[\sigma_K]=c_{top}({\Sym}^kU^*)
\pi_{l*}(c_{top}(\pi_l^*{\Sym}^{d}U^*/(\pi_l^*{\Sym}^kU^*\otimes T_l))
c_{top}(\pi_l^*{\Sym}^lU^*/T_l)) .\tag 3.2
$$
The corresponding formula holds for the class of $\sigma_L$,
that is,
$$
[\sigma_L]=c_{top}({\Sym}^lU^*)
\pi_{k*}(c_{top}(\pi_k^*{\Sym}^{d}U^*/(\pi_k^*{\Sym}^lU^*\otimes T_k))
c_{top}(\pi_k^*{\Sym}^kU^*/T_k)) .\tag 3.3
$$
\endproclaim{}

\demo{Proof of Theorem 3.1}
Recall that
$$
\eqalign{\sigma_K=&\{\Bbb P^r \in G(r+1,n+1)\ |
\  \Bbb P^r \subset K,\ \ \Bbb P^r \cap L \subset D\}\cr
=& \{\Bbb P^r \in G(r+1,n+1) \ | \  \Bbb P^r \subset K\}
\cap \{\Bbb P^r \in G(r+1,n+1) \ | \  \Bbb P^r \cap L \subset D\}.\cr}
$$
Since $K$ is a generic hypersurface of degree $k$,
$$
[\sigma_K]=c_{top}({\Sym}^kU^*)[B(D,L)],
$$
where $[B(D,L)]$ is the class of
$$
B(D,L)=\{\Bbb P^r \in G(r+1,n+1) \ | \  \Bbb P^r \cap L \subset D\}.
$$
To see $[B(D,L)]$,
we will first rewrite $\Bbb P({\Sym}^lU^*)$ as
$$
\Bbb P({\Sym}^lU^*)=\{(X_l, \Bbb P^r)\ | \ X_l \subset \Bbb P^r\}.
$$
Notice that,
for $l > 1$,
$\Bbb P^r$ in $(X_l,\Bbb P^r)$ is uniquely determined by $X_l$
and we will call this $\Bbb P^r$ the supporting plane of $X_l$.
Therefore,
geometrically,
we can identify $\Bbb P({\Sym}^lU^*)$ as the space
of degree $l$ subvarieties $X_l$ of dimension $(r-1)$ in
$\Bbb P^n$ such that
$X_l$ is contained in some $\Bbb P^r$.
In the case of $\ l=1$,
$\ \Bbb P({\Sym}^lU^*)$ can be identified as the space of
marked $\Bbb P^{r-1}$'s in $\Bbb P^n$ with the mark
to be a $\Bbb P^r$ assigned to
each $\Bbb P^{r-1}$.
Such a $\Bbb P^r$ is taken from all $\Bbb P^r$'s
containing the given $X_1=\Bbb P^{r-1}$
and we will also call it the supporting plane of $X_1$.
With such an understanding,
we will from now on denote an element of $\Bbb P({\Sym}^lU^*)$
just by $X_l$.
Now,
consider the commutative diagram
$$
\CD
B^*(D,L) @>f>> \Bbb P({\Sym}^lU^*) \\
@VViV @VV\pi_lV\\
B(D,L) @>j>> G(r+1,n+1)
\endCD{}
$$
where
$$
B^*(D,L)=B(D,L)-\{\Bbb P^r\ | \ \Bbb P^r \subset L\},
$$
$i$ and $j$ are inclusions, and $f$ is defined by
$$
f(\Bbb P^r)=\Bbb P^r \cap L.
$$
{}From the definition,
we see that $f$ is $1$-to-$1$.
Moreover,
$f$ can be considered as
a rational map on $B(D,L)$.
Since
$$
\pi_l\circ f |_{B^*(D,L)}= id,
$$
we therefore have that
$$
\aligned
[B(D,L)]=
&[\overline{B^*(D,L)}]\\
=&[\overline{\pi_l \circ f(B^*(D,L))}]\\
= & [\pi_l (\overline{f(B^*(D,L))})]\\
=&\pi_{l*} ([\overline{f(B^*(D,L))}]).
\endaligned{}
\tag 3.4
$$

The following lemma is needed to complete our argument
and its proof will be given later.

\proclaim{Lemma 3.5}
$$
\overline{f(B^*(D,L))}=
\{X_l\ | \ X_l \subset D\}
\cap \{X_l\ | \ X_l \subset L\}.
$$
\endproclaim{}

To see how to use Lemma 3.5 to derive formula (3.2),
consider the subbundle
$$
\pi_l^*{\Sym}^kU^*\otimes T_l \ \hookrightarrow
\ \pi_l^*{\Sym}^{d}U^*
$$
on $\Bbb P({\Sym}^lU^*)$.
(Recall that $d=k+l$.)
This defines a quotient bundle $Q_d$ on $\Bbb P({\Sym}^lU^*)$ which fits
into the exact sequence
$$
0 \to \pi_l^*{\Sym}^kU^*\otimes T_l
\to \pi_l^*{\Sym}^{d}U^* \to Q_d \to 0.
$$
Hypersurface $D$ induces a section of ${\Sym}^{d}U^*$
and hence a section $s_D$ of $Q_d$.
The zero-scheme of $s_D$ consists exactly of $X_l$'s
which are contained in $D$.
In fact,
at each $X_l$ in $\Bbb P({\Sym}^lU^*)$,
the value of the section induced by $D$ in
$\pi_l^*{\Sym}^{d}U^*$ is given by the
intersection of $D$
with the supporting plane of $X_l$.
It lies in $\pi_l^*{\Sym}^kU^*\otimes T_l$ if and only if this intersection
is the union of $X_l$ and
some degree $k$ hypersurface in the supporting plane.
On the other hand,
it is easy to see that $Q_d$ is generated by sections
induced from hypersurfaces of degree $d$ in $\Bbb P^n$.
Since $D$ is generic,
$s_D$ is therefore a regular section of $Q_d$.
Hence,
by the localized top Chern class theorem,
$$
[\{X_l\ | \ X_l \subset D\}]=[Z(s_D)]
=c_{top}(Q_d)=c_{top}(\pi_l^*{\Sym}^{d}U^*/(\pi_l^*{\Sym}^kU^*\otimes T_l)).
$$
In a similar way,
we have that
$$
[\{X_l\ | \ X_l \subset L\}]=
c_{top}(\pi_l^*{\Sym}^lU^*/T_l).
$$
Lemma 3.5 and formula (3.4) together now yield the formula for $[B(D,L)]$
and hence formula (3.2) for $[\sigma_K]$, that is,
$$
\aligned
[\sigma_K]=&c_{top}({\Sym}^kU^*)[B(D,L)]\\
=&c_{top}({\Sym}^kU^*)
\pi_{l*}(c_{top}(\pi_l^*{\Sym}^{d}U^*/(\pi_l^*{\Sym}^kU^*\otimes T_l))
c_{top}(\pi_l^*{\Sym}^lU^*/T_l)).
\endaligned{}
$$

To complete our proof,
we still need to show Lemma 3.5 that

$$
\overline {f(B^*(D,L))}=
\{X_l\ | \ X_l \subset D\}
\cap \{X_l\ | \ X_l \subset L\}.
\tag 3.6
$$
To simplify the notation,
we will call the right-hand side of (3.6) $B$.
Recall that
$$
Q_d=\pi_l^*{\Sym}^{k+l}U^*/(\pi_l^*{\Sym}^kU^*\otimes T_l).
$$
Let also
$$
Q_l=\pi_l^*{\Sym}^lU^*/T_l.
$$
As we have shown,
$B$ is the intersection of the zero-scheme of $s_D$ and
the zero-scheme of $s_L$,
where $s_D$ and $s_L$ are sections of
$Q_d$ and $Q_l$
induced by $D$ and $L$,
respectively.
Moreover,
those two quotient bundles
$Q_d$
and
$Q_l$
are generated by sections
induced by hypersurfaces of degree $d$ and degree $l$ in $\Bbb P^n$,
respectively.
Since $D$ and $L$ are generic,
we see that $B$ is a smooth subvariety of pure dimension in $\Bbb
P({\Sym}^lU^*)$.
Let
$$
B^*=\{X_l \in B\ | \ \text{the supporting plane of $X_l$
is not contained in $L$}\}.
$$
It is easy to see from the definition of $f$ that
$$
f(B^*(D,L))=B^*.
$$
Therefore,
to show (3.6),
it is enough to show that $B-B^*$ is proper in $B$.
To see this,
notice that
$$
\aligned
B-B^*=&\{X_l\ | \ X_l \subset D\}
\cap \{X_l\ | \ X_l \subset \Bbb P^r \subset L\}\\
=&\{X_l\ | \ X_l \subset D\}
\cap (\bigcup_{\Bbb P^r \subset L} \pi_l^{-1}(\Bbb P^r)).
\endaligned{}
$$
Since $D$ is generic,
it is enough for us to show that
$$
\dim(\bigcup_{\Bbb P^r \subset L} \pi_l^{-1}(\Bbb P^r)) <
\dim(\{X_l\ | X_l \subset L\}).
$$
This is easy to see.
In fact,
we have either
$$
\bigcup_{\Bbb P^r \subset L} \pi^{-1}(\Bbb P^r)= \emptyset
$$
or
$$
\aligned
\dim(\bigcup_{\Bbb P^r \subset L} \pi^{-1}(\Bbb P^r))=&\rank({\Sym}^lU^*)-1+
(\dim(G(r+1,n+1))-\rank({\Sym}^lU^*))\\
=&\dim(G(r+1,n+1))-1.
\endaligned{}
$$
On the other hand,
$$
\aligned
\dim(\{X_l\ | \ X_l \subset L\})=&\dim(\Bbb P({\Sym}^lU^*))-\rank(Q_l)\\
=&\dim(G(r+1,n+1))+\rank({\Sym}^lU^*)-1 - (\rank({\Sym}^lU^*)-1)\\
=&\dim(G(r+1,n+1)).
\endaligned{}
$$
Therefore,
$B-B^*$ is either empty or a divisor in $B$.
This completes our proof of Lemma 3.5
and hence our proof of formula (3.2).
By interchanging $l$ and $k$,
the same argument will yield formula (3.3) for the class of $[\sigma_L]$.
\qed
\enddemo{}

Theorem 3.1 gives us a tool to compute
the classes of $[\sigma_K]$ and $[\sigma_L]$
using standard techniques in intersection theory.
We have made some computations using a Maple package
written by Katz and Str{\o}mme.
Computations of examples have
played important roles in discovering Theorem 3.1.
In fact,
many of them were made before Theorem 3.1 was proved.
A part of those examples along with other things will be given in Section 5.

\smallpagebreak
\heading 4. Summary of main results obtained so far
\endheading{}
\smallpagebreak

We will summarize the main results obtained so far in this section.
For more details on notation,
please refer to corresponding theorems and formulas as indicated.

\proclaim{Summary/Theorem 4.1}
For a generic degeneration
$$
\Cal D =\{X_s=KL+sD+\cdots\}_{s \in \Delta},
$$
there are no higher order obstructions to deform $\Bbb P^r$
with $X_0$.
The subscheme $F_{\text {\, lim}}$ of limiting $\Bbb P^r$'s is
hence determined by the infinitesimal data.
In particular,
$$
F_{\text {\, lim}}=\sigma_K\cup \sigma_L ,
$$
where,
as in (2.5),
$$
\sigma_K = \{\Bbb P^r \in G(r+1,n+1)|
\  \Bbb P^r \subset K, \ \Bbb P^r \cap L \subset D\}
$$
and
$$
\sigma_L = \{\Bbb P^r \in G(r+1,n+1)|
\  \Bbb P^r \subset L, \ \Bbb P^r \cap K \subset D\}
$$
consist of limiting $\Bbb P^r$'s in $K$ and $L$,
respectively.
Furthermore,
their classes can be computed by formulas (3.2) and (3.3)
$$
[\sigma_K]=c_{top}({\Sym}^kU^*)
\pi_{l*}(c_{top}(\pi_l^*{\Sym}^{d}U^*/(\pi_l^*{\Sym}^kU^*\otimes T_l))
c_{top}(\pi_l^*{\Sym}^lU^*/T_l))
$$
and
$$
[\sigma_L]=c_{top}({\Sym}^lU^*)
\pi_{k*}(c_{top}(\pi_k^*{\Sym}^{d}U^*/(\pi_k^*{\Sym}^lU^*\otimes T_k))
c_{top}(\pi_k^*{\Sym}^kU^*/T_k)).
$$
\endproclaim{}

By using the setup in Section 1,
the above facts can now be applied to the refined intersection product
defined there.
Notice that the following statement is purely pointwise on fixed $X_0$
and no deformation is involved.

\proclaim{Summary/Theorem 4.1'}
Let ${R}_{X_0}$ be the refined intersection product
in $A_m(F_{X_0})$ as defined in Section 1.
The contributions (equivalences) of $F_K$ and $F_L$ to ${R}_{X_0}$
are
$$
{R}_{X_0}^{F_K}=c_{top}({\Sym}^kU^*)
\pi_{l*}(c_{top}(\pi_l^*{\Sym}^{d}U^*/(\pi_l^*{\Sym}^kU^*\otimes T_l))
c_{top}(\pi_l^*{\Sym}^lU^*/T_l))
$$
and
$$
{R}_{X_0}^{F_L}=c_{top}({\Sym}^lU^*)
\pi_{k*}(c_{top}(\pi_k^*{\Sym}^{d}U^*/(\pi_k^*{\Sym}^lU^*\otimes T_k))
c_{top}(\pi_k^*{\Sym}^kU^*/T_k)),
$$
respectively.
In particular,
we have the canonical decomposition of ${R}_{X_0}$ as follows
$$
\aligned
{R}_{X_0}=&c_{top}({\Sym}^kU^*)
\pi_{l*}(c_{top}(\pi_l^*{\Sym}^{d}U^*/(\pi_l^*{\Sym}^kU^*\otimes T_l))
c_{top}(\pi_l^*{\Sym}^lU^*/T_l))\\
&+c_{top}({\Sym}^lU^*)
\pi_{k*}(c_{top}(\pi_k^*{\Sym}^{d}U^*/(\pi_k^*{\Sym}^lU^*\otimes T_k))
c_{top}(\pi_k^*{\Sym}^kU^*/T_k)).
\endaligned{}
\tag 4.2
$$
\endproclaim{}

\demo{Proof}
Except in the case of $r=1$,
it may be possible that ${R}_{X_0}$
will have $F_K\cap F_L$ as a distinguished
variety.
However,
from the definitions of $\sigma_K$ and $\sigma_L$,
we see that the subscheme of limiting
$\Bbb P^r$'s contained in $K\cap L$ must have a dimension strictly
less than $m$.
Therefore,
the equivalence of $F_K\cap F_L$ would be zero in any case.
This implies that the sum of $[\sigma_K]$ and $[\sigma_L]$
is indeed a decomposition of ${R}_{X_0}$.
\qed
\enddemo{}

\example{Question}
There is one more decomposition of ${R}_{X_0}$ that is worth
mentioning.
If $F_K$ was a connected component of $F_{X_0}$,
then the contribution of $F_K$ to ${R}_{X_0}$
would be
$$
{R}_{X_0}^{F_K}=\{c({\Sym}^dU^*)\cup s(F_K,G)\}_m.
$$
This class is very easy to compute.
In general,
we can still decompose ${R}_{X_0}$ as the sum of
the class above and the so-called residual intersection class $\Bbb R$.
This is the celebrated residual intersection formula of
Fulton-Laksov [F, Chapter 9] [FL] [L].
We can also compute
the similar residual decomposition of ${R}_{X_0}$ for $F_L$.
What are geometric meanings
for such decompositions?
We have computed some examples but could not give satisfactory answers.
\endexample{}

\smallpagebreak
\heading 5. Identities in characteristic classes and examples
\endheading{}
\smallpagebreak

In this section,
we will derive some identities
in characteristic classes from the results in previous sections
and then compute a few examples with our formulas.
{}From the splitting principle,
those identities can be translated into identities in
symmetric polynomials of $r+1$ variables and hence should
hold for any vector bundle of rank $r+1$ over a manifold.
It will be interesting to see if the identities can be verified
by a more direct approach.

Since the image of ${R}_{X_0}$ in the Chow ring of $G(r+1,n+1)$
is always equal to the top Chern class
of ${\Sym}^dU^*$,
the results in the last section
immediately imply the following family of identities.

\proclaim{Corollary 5.1}
For any $d$, $k$, and $l$ with $d=k+l$,
the following identity holds on the Chow ring of $G(r+1,n+1)$.
$$
\aligned
&c_{top}({\Sym}^dU^*)\\
=&c_{top}({\Sym}^kU^*)
\pi_{l*}(c_{top}(\pi_l^*{\Sym}^{d}U^*/(\pi_l^*{\Sym}^kU^*\otimes T_l))
c_{top}(\pi_l^*{\Sym}^lU^*/T_l))\\
+&c_{top}({\Sym}^lU^*)
\pi_{k*}(c_{top}(\pi_k^*{\Sym}^{d}U^*/(\pi_k^*{\Sym}^lU^*\otimes T_k))
c_{top}(\pi_k^*{\Sym}^kU^*/T_k)).
\endaligned{}
\tag 5.2
$$
\endproclaim{}

Before giving examples,
we will first rewrite formula (3.2) and formula (3.3)
more explicitly in terms of elements in the Chow ring of $G(r+1,n+1)$.
One of the reasons for doing this is because
it will be useful in actual computations.
Since the Chow ring of the projective bundle $\Bbb P({\Sym}^lU^*)$
is bigger and more complicated than the Chow ring of the underlying space
$G(r+1,n+1)$,
in many cases,
a direct computation on $G(r+1,n+1)$ will be faster than
computing first on $\Bbb P({\Sym}^lU^*)$ and then pushing it down.
More importantly,
by computing on $G(r+1,n+1)$ only,
we eliminate the need to create the projective bundle
$\Bbb P({\Sym}^lU^*)$,
which may take too much time and space on a computer,
especially if the rank of ${\Sym}^lU^*$ is big.

For a given vector bundle $E$,
we will use $c_i(E)$ to denote the $i^{th}$ Chern class of $E$
and $s_i(E)$ to denote the $i^{th}$ Segre class of $E$.
Recall also that we identify an element in $A^*(G(r+1,n+1))$ with its dual in
$A_*(G(r+1,n+1))$.

\proclaim{Proposition 5.3}
Formula (3.2) can be rewritten directly in terms of elements in
the Chow ring of $G(r+1,n+1)$ as
$$
\aligned
[\sigma_K]=&c_{r_k}({\Sym}^kU^*)
\sum_{i=0}^{r_d-r_k}
\sum_{j=0}^{r_l-1}
\sum_{h=0}^{r_d-r_k-i-j}
{r_d-1-i \choose r_k-1+h}
c_{i}({\Sym}^{d}U^*)\\
&\times c_j({\Sym}^lU^*)s_{h}({\Sym}^kU^*)s_{r_d-r_k-h-i-j}({\Sym}^lU^*) ,
\endaligned{}
\tag 5.4
$$
where
$$
\left\{\aligned
r_d=&\rank({\Sym}^dU^*)={r+d \choose r},\\
r_k=&\rank({\Sym}^kU^*)={r+k \choose r},\\
r_l=&\rank({\Sym}^lU^*)={r+l \choose r}.
\endaligned{}
\right .
$$
The corresponding formula holds for $[\sigma_L]$ if we interchange
$k$ and $l$.
\endproclaim{}

\example{Remark}
Since $s_i(E)=0$ for all $i<0$,
we may replace the range of $j$ by
$$
0 \leq j \leq \min\{r_l-1, r_d-r_k-i\}.
$$
\endexample{}

\demo{Proof}
This is a straightforward computation with standard
calculations on the Chern classes and the Segre classes.
All the formulas used here can be found in Fulton's book [F].
{}From formula (3.2),
$$
\aligned
[\sigma_K]=&c_{top}({\Sym}^kU^*)
\pi_{l*}(c_{top}(\pi_l^*{\Sym}^{d}U^*/(\pi_l^*{\Sym}^kU^*\otimes T_l))
c_{top}(\pi_l^*{\Sym}^lU^*/T_l))\\
=& c_{r_k}({\Sym}^kU^*)
\pi_{l*}(\sum_{i=0}^{r_d-r_k}
c_{i}(\pi_l^*{\Sym}^{d}U^*)s_{r_d-r_k-i}(\pi_l^*{\Sym}^kU^*\otimes T_l)\\
&\times \sum_{j=0}^{r_l-1}
c_{j}(\pi_l^*{\Sym}^lU^*)s_{r_l-1-j}(T_l)) \\
=& c_{r_k}({\Sym}^kU^*)
\pi_{l*}(\sum_{i=0}^{r_d-r_k}c_{i}(\pi_l^*{\Sym}^{d}U^*)
\sum_{h=0}^{r_d-r_k-i}{r_k-1+r_d-r_k-i\choose r_k-1+h}\\
&\times s_{h}(\pi_l^*{\Sym}^kU^*)c_1(-T_l)^{r_d-r_k-i-h}
\sum_{j=0}^{r_l-1}
c_{j}(\pi_l^*{\Sym}^lU^*)c_1(-T_l)^{r_l-1-j}) \\
=& c_{r_k}({\Sym}^kU^*)
\sum_{i=0}^{r_d-r_k}
\sum_{j=0}^{r_l-1}
\sum_{h=0}^{r_d-r_k-i}{r_d-1-i\choose r_k-1+h}
c_{i}({\Sym}^{d}U^*)\\
&\times
c_{j}({\Sym}^lU^*)s_{h}({\Sym}^kU^*)\pi_{l*}(c_1(\Cal
O(1))^{r_d-r_k+r_l-1-i-j-h})\\
=& c_{r_k}({\Sym}^kU^*)
\sum_{i=0}^{r_d-r_k}
\sum_{j=0}^{r_l-1}
\sum_{h=0}^{r_d-r_k-i-j}{r_d-1-i\choose r_k-1+h}
c_{i}({\Sym}^{d}U^*)\\
&\times c_{j}({\Sym}^lU^*)s_{h}({\Sym}^kU^*)s_{r_d-r_k-i-j-h}({\Sym}^lU^*).
\endaligned{}
$$
\qed
\enddemo{}

In the case of $r=1$,
we have
$$
r_d=d+1,\quad r_k=k+1,\quad r_l=l+1.
$$

\proclaim{Corollary 5.5}
For the case of lines,
$$
\aligned
[\sigma_K]=&c_{k+1}({\Sym}^kU^*)
\sum_{i=0}^{l}
\sum_{j=0}^{l-i}
\sum_{h=0}^{l-i-j}
{d-i \choose k+h}
c_{i}({\Sym}^{d}U^*)c_j({\Sym}^lU^*)\\
&\times s_{h}({\Sym}^kU^*)s_{l-h-i-j}({\Sym}^lU^*).
\endaligned{}
\tag 5.6
$$
The corresponding formula holds for $[\sigma_L]$ if we interchange
$k$ and $l$.
\endproclaim{}

Corollary 5.1 can now be restated as well.

\proclaim{Corollary 5.7}
For any set of non-negative integers
$r$, $n$, $d$, $k$, and $l$ with $r < n$ and $d=k+l$,
the following identity holds on the Chow ring of $G(r+1,n+1)$.
$$
\aligned
c_{top}({\Sym}^dU^*)=&
c_{r_k}({\Sym}^kU^*)
\sum_{i=0}^{r_d-r_k}
\sum_{j=0}^{r_l-1}
\sum_{h=0}^{r_d-r_k-i-j}
{r_d-1-i \choose r_k-1+h}
c_{i}({\Sym}^{d}U^*)\\
&\times c_j({\Sym}^lU^*)s_{h}({\Sym}^kU^*)s_{r_d-r_k-h-i-j}({\Sym}^lU^*)\\
&+c_{r_l}({\Sym}^lU^*)
\sum_{i=0}^{r_d-r_l}
\sum_{j=0}^{r_k-1}
\sum_{h=0}^{r_d-r_l-i-j}
{r_d-1-i \choose r_l-1+h}
c_{i}({\Sym}^{d}U^*)\\
&\times c_j({\Sym}^kU^*)s_{h}({\Sym}^lU^*)s_{r_d-r_l-h-i-j}({\Sym}^kU^*) ,
\endaligned{}
\tag 5.8
$$
where
$$
\left\{\aligned
r_d=&\rank({\Sym}^dU^*)={r+d \choose r},\\
r_k=&\rank({\Sym}^kU^*)={r+k \choose r},\\
r_l=&\rank({\Sym}^lU^*)={r+l \choose r}.
\endaligned{}
\right .
$$
Moreover,
this gives the canonical decomposition of ${R}_{X_0}$.
\endproclaim{}

\example{Remark}
We think identities (5.2) and (5.8) should hold for
any vector bundle of rank $r+1$
over a general manifold.
In fact,
no relations in Chow rings should be needed and
they should hold as identities of symmetric polynomials
as we have seen in [W] for the case
of $r=1$.
However,
we do not know whether such identities would have further significant meanings
in general.
It will be very interesting if one can find other applications
of such identities in different settings.
\endexample{}

We will now compute some examples using either Theorem 3.1
or Proposition 5.3.
Most of our computations are done by using the schubert package [KS]
on Maple.
It is very easy to write a schubert code based on formula (3.2)
and formula (3.3).
On the other hand,
although formulas given in Proposition 5.3 are more
complicated,
their computations in many cases are much faster.
To keep track of the degrees of the hypersurfaces in our degeneration,
we will use ${R}(d,k)$ for
$[\sigma_K]$ and ${R}(d,l)$ for $[\sigma_L]$
in the following examples.
In theory,
one can compute any ${R}(d,k)$ and ${R}(d,l)$ in $A^*(G(r+1,n+1))$.
In practice,
however,
one will find out quickly that $d$ or $r$ cannot be set too large
even with the help of powerful computers.

\proclaim {Example 1} Degenerations of hypersurfaces and
their lines ($r=1$).
\endproclaim{}
In this case,
we got the same results as given in [W] and
one can find detailed examples,
such as degenerations of cubic surfaces, or quintic threefolds,
there.
Since the formulas used here are not the same
as given in [W],
it is a reassurance to see that we do get the same results.
Comparing Corollary 5.5 of this paper and Theorem 1.1 in [W],
we see that the formulas here are more complicated
and have many extra terms.
Those extra terms must actually cancel out each other.
However,
it does not appear to be trivial to prove directly that the formulas
given in Corollary 5.5 are actually equal to
the corresponding formulas in [W].
Once again,
we think that they should be equal to each other as
polynomials and hence the equations should hold for general vector bundles
of rank 2.

\proclaim {Example 2} Degenerations of
a generic
quartic hypersurface in $\Bbb P^7$ and its $\Bbb P^2$'s.
\endproclaim{}
See also [K3] for a computation by geometric methods.
It is that computation which leads to the discovery of our formulas.

There are $3,297,280$
$\Bbb P^2$'s on a generic
quartic hypersurface in $\Bbb P^7$.
This is given by the degree of $c_{top}({\Sym}^4U^*)$
on $A^*(G(3,8))$.
There are two degenerations
and computations using our representations show that
\roster
\itemitem{} Case 1, $K$ is a cubic and $L$ is a plane:
$\quad {R}(4,3)=483,840$,
${R}(4,1)=2,813,440$.
\itemitem{} Case 2, $K$ and $L$  are quadrics:
$\quad \quad \quad \quad \ {R}(4,2)=1,648,640$.
\endroster{}
Notice that in each of the cases above we do have that
$$
{R}(4,k)+{R}(4,l)=3,297,280=c_{top}({\Sym}^7U^*)
$$
as expected.

\proclaim {Example 3} Degenerations of
a generic quadric in $\Bbb P^n$, $n\geq 3$, and its $\Bbb P^2$'s
\endproclaim{}

To eliminate doubtful feelings one might have
towards computers,
this example has been double checked by hand and details are
given below.

Let
$$
c_1(U^*)=x,\quad
c_2(U^*)=y,\quad
c_3(U^*)=z.
$$
It is easy to compute that
$$
\left\{
\aligned
c_1({\Sym}^2U^*)=&4x,\\
c_2({\Sym}^2U^*)=&5(x^2+y),\\
c_3({\Sym}^2U^*)=&2x^3+11xy+7z,\\
\vdots \quad & \quad \vdots \\
c_6({\Sym}^2U^*)=&8z(xy-z).
\endaligned{}
\right .
\tag 5.9
$$
Therefore,
the class of $\Bbb P^2$'s on a generic
quadric is equal to
$$
8xyz-8z^2.
$$
On the other hand,
to compute the class of limiting $\Bbb P^2$'s
when the quadric degenerates into two $(n-1)$-planes $K$ and $L$,
we set
$$
r=d=2, \quad l=k=1
$$
in (5.4) of Proposition 5.3.
This gives us
$$
{R}(2,1)=c_{3}(U^*)
\sum_{i=0}^{3}
\sum_{j=0}^{2}
\sum_{h=0}^{3-i-j}
{5-i \choose 2+h}
c_{i}({\Sym}^{2}U^*)c_j(U^*)s_{h}(U^*)s_{3-h-i-j}(U^*).
\tag 5.10
$$
It is easy to see that
$$
\left\{
\aligned
s_1({\Sym}^2U^*)=&-x,\\
s_2({\Sym}^2U^*)=&x^2-y,\\
s_3({\Sym}^2U^*)=&-x^3+2xy-z,
\endaligned{}
\right .
\tag 5.11
$$
Substituting (5.9) and (5.11) into (5.10) and simplifying,
we finally have
$$
{R}(2,1)=4xyz-4z^2
$$
as expected.

We have purposefully done our computations without using the Schubert
calculus.
Therefore,
this also verifies the generalization of identity (5.8) to any
vector bundle of rank $3$ in the case of $d=2$.

\example{Remark}
It is well-known that there is no $\Bbb P^2$ in a generic quadric in $\Bbb P^n$
if $n < 5$.
This can be seen by using the Schubert calculus.
On $A_*(G(3,n+1))$,
$$
c_6({\Sym}^2U^*)=8z(xy-z)=8\sigma_{3,2,1}.
$$
(In fact,
$8\sigma_{3,2,1}$
is equal to the class of $\Bbb P^2$'s in any smooth quadric.
See [GH, Chapter 6] for computations using geometric methods.)
Hence,
this class is zero if $n<5$.
\endexample{}

Our last example shows something interesting.

\proclaim {Example 4} Degenerations of a generic
cubic hypersurface in $\Bbb P^8$ and its $\Bbb P^3$'s.
\endproclaim{}

By computing the degree of $c_{top}({\Sym}^3U^*)$
on $A^*(G(4,9))$,
we see that there are $321,489$
$\Bbb P^3$'s on a generic
cubic hypersurface in $\Bbb P^8$.
When the cubic degenerates into the union of a quadric and a plane,
computations from our formulas yield that
$$
{R}(3,2)=0, \quad {R}(3,1)=321,489.
$$
The result is a surprise since it says that all limiting $\Bbb P^3$'s
are in the plane and none of them is in the quadric.
In terms of intersection theory,
that means the equivalence of $F_K$ for ${R}_{X_0}$ is zero.
Notice that ${\Sym}^dU^*$ is a bundle generated by sections
so all equivalences must be non-negative classes.
Our example gives the extreme case
and implies that ${\Sym}^dU^*$ (and hence $U^*$) is not a positive bundle.

\smallpagebreak
\heading Appendix: Schubert code
\endheading{}
\smallpagebreak

In this appendix,
we will give schubert code for Example 3 and Example 4 of Section 5.
Formula (4.4) of Proposition 5.3 will be applied for Example 3 and
Formula (3.2) and formula (3.3) of Theorem 3.1 will be used
for Example 4.
By changing corresponding parameters $r$, $n$, $d$, $k$, and $l$,
those programs can be converted easily to compute other examples.
Schubert or Maple commands that a user is required to type
during a Maple session will be indicated by the line beginning
with the Maple prompt ``$>$''.
Comments will be indicated by the line beginning with ``\#''.
In many places,
we have used the command ``latex'' to make tex-format output.
One should leave them out if regular output is desired.
We want to thank Katz and Str{\o}mme
for making this wonderful Maple package available.
For those who may be interested in getting the current version of schubert,
we want to point out that it is available
by using ftp to ftp.math.okstate.edu with
user$=$anonymous,
password$=<$your email address$>$,
cd pub/schubert,
then get the files there.
One can find instructions on how to install the package
from the manual.

\example{A note of caution}
Schubert uses notation from the school of Grothendieck.
Since notation in this paper primarily follows from Fulton's book,
we have to use the dual to covert between them in many places.
This may cause some confusion.
It could easily lead to wrong answers if one is not careful about this.
\endexample{}

\proclaim {Schubert code for Example 3} Degenerations of
a generic quadric in $\Bbb P^n$, $n=4$, and its $\Bbb P^2$'s.
Results will be same for all $n>4$ except for
the last computation.
\endproclaim{}

\# Construct $G(3,5)=Gc$ with $U^* = Qc$ and $c_i(U^*) = c_i$.

$>$ latex(grass(3,5,c,mon));

$$
currentvariety\_ \ is \ Gc, \ \ DIM \ is \ 6
$$

\# Construct ${\Sym}^2U^*=s2$.

$>$ s2:=symm(2,Qc):

\# The class of $\Bbb P^2$'s
in a generic quadric is equal to the top chern class
of ${\Sym}^2U^*$.

$>$ cl\_P2\_in\_X2:=chern(rank(s2),s2):

$>$ latex($''$);

$$
cl\_P2\_in\_X2\ = \ 8\,{c_3}\,{c_1}\,{c_2}-8\,{c_3}^{2}
$$

\# Use formula (5.10) to compute the class of limiting $\Bbb P^2$'s
in a hyperplane which should

\# be equal to one half of the class above.

$>$ f:=0:

$>$ for i from 0 to 3 do

$>$ for j from 0 to min(2,3-i) do

$>$ for h from 0 to 3-i-j do

$>$ f:=f+

$>$ binomial(5-i,2+h)$*$chern(i,s2)$*$chern(j,Qc)
$*$segre(h,dual(Qc))$*$segre(3-h-i-j,dual(Qc))

$>$ od od od:

$>$ cl\_limiting\_P2\_in\_X1:=expand(chern(3,Qc)$*$f):

$>$ latex($''$);

$$
cl\_limiting\_P2\_in\_X1 \ = \ 4\,{c_3}\,{c_1}\,{c_2}-4\,{c_3}^{2}
$$

\# In our case of $n=4$, $\dim(G(3,5))=\rank({\Sym}^2U^*)$,
so the class above gives a number.

$>$ number:=integral(Gc,cl\_limiting\_P2\_in\_X1):

$>$ latex($''$);

$$
number\ = \ 0
$$

$>$ quit;

\proclaim {Schubert code for Example 4}
Degenerations of a generic cubic hypersurface in $\Bbb P^8$
and its $\Bbb P^3$'s.
\endproclaim{}

\# Construct $G(4,9)=Gc$ with $U^* = Qc$ and $c_i(U^*) = c_i$.

$>$ latex(grass(4,9,c,mon));

$$
currentvariety\_ \ is \ Gc, \ DIM \ is \ 20
$$

\# Construct ${\Sym}^3U^*=s3$.

$>$ s3:=symm(3,Qc):

\# The class of $\Bbb P^3$'s in a generic cubic hypersurface
is equal to the top chern class

\# of ${\Sym}^3U^*$.

$>$ cl\_P3\_in\_X3:=chern(20,s3):

$>$ latex($''$);
$$
\aligned
cl\_P3\_in\_X3 \, = \,
&-1296\,{c_4}^{2}{c_3}^{3}{c_1}^{3}
-2592\,{c_4}\,{c_3}^{3}{c_1}\,{c_2}^{3}
+17496\,{c_4}^{2}{c_3}^{4}\\
&+1296\,{c_4}^{3}{c_2}^{4}
+15552\,{c_4}^{2}{c_3}^{2}{c_2}\,{c_1}^{4}
+2592\,{c_4}\,{c_3}^{2}{c_2}^{4}{c_1}^{2}\\
&+50625\,{c_4}^{5}
+2592\,{c_4}^{2}{c_3}^{2}{c_2}^{3}
-14580\,{c_4}^{3}{c_3}^{2}{c_2}\\
&+36045\,{c_4}^{4}{c_2}\,{c_1}^{2}
+17496\,{c_4}^{3}{c_2}\,{c_1}^{6}
+34425\,{c_4}^{4}{c_1}\,{c_3}\\
&+1296\,{c_4}^{2}{c_2}^{5}{c_1}^{2}
-17496\,{c_4}^{3}{c_3}\,{c_1}^{5}
-81162\,{c_4}^{3}{c_3}^{2}{c_1}^{2}\\
&-13608\,{c_4}^{3}{c_2}^{3}{c_1}^{2}
+2592\,{c_4}^{2}{c_2}^{4}{c_1}^{4}
-14580\,{c_4}^{3}{c_2}^{2}{c_1}^{4}\\
&-17496\,{c_4}\,{c_1}\,{c_3}^{5}
-2592\,{c_4}\,{c_3}^{4}{c_1}^{4}
+87966\,{c_4}^{3}{c_3}\,{c_2}\,{c_1}^{3}\\
&+3888\,{c_4}^{3}{c_3}\,{c_1}\,{c_2}^{2}
-16200\,{c_4}^{4}{c_2}^{2}
+17496\,{c_4}^{4}{c_1}^{4}\\
&+2916\,{c_4}^{2}{c_3}^{2}{c_2}^{2}{c_1}^{2}
-1296\,{c_4}^{2}{c_3}\,{c_1}\,{c_2}^{4}
-11664\,{c_4}^{2}{c_3}\,{c_2}^{2}{c_1}^{5}\\
&+2916\,{c_4}^{2}{c_3}^{3}{c_1}\,{c_2}
-14904\,{c_4}^{2}{c_3}\,{c_2}^{3}{c_1}^{3}
+5184\,{c_4}\,{c_3}^{2}{c_2}^{3}{c_1}^{4}\\
&+29160\,{c_4}\,{c_3}^{4}{c_2}\,{c_1}^{2}
+2592\,{c_4}\,{c_3}^{3}{c_2}\,{c_1}^{5}
-16848\,{c_4}\,{c_3}^{3}{c_2}^{2}{c_1}^{3}
\endaligned{}
$$

\# Since $\dim(G(4,9))=\rank({\Sym}^3U^*)$,
the class above gives the number of $\Bbb P^3$'s in

\# a generic cubic hypersurface in $\Bbb P^8$.

$>$ total\_number:=integral(Gc,cl\_P3\_in\_X3):

$>$ latex($''$);

$$
total\_number\ = \ 321489
$$

\# Construct ${\Sym}^2U^*=s2$.

$>$ s2:=symm(2,Qc):

\# Construct the projective bundle $\Bbb P(U^*)$ with
the projection map $\pi_1 = P1$ and the

\# tautological bundle $T_1 = z1$.

$>$ Proj(P1,dual(Qc),z1):

\# Use formula (3.2) to compute the class and the number
of limiting $\Bbb P^3$'s in a generic

\# quadric.

$>$ cl\_limiting\_P3\_in\_X2:=

$>$ expand(chern(10,s2)$*$
lowerstar(P1,chern(10,s3-s2{\&}$*$o(-z1))$*$chern(3,Qc-o(-z1)))):

$>$ latex($''$);

$$
\aligned
cl\_limiting\_P3\_in\_X2\, = \,
&-12528\,{c_4}^{2}{c_3}^{3}{c_1}^{3}
-1728\,{c_4}\,{c_3}^{3}{c_1}\,{c_2}^{3}
+11664\,{c_4}^{2}{c_3}^{4}\\
&+17280\,{c_4}^{2}{c_3}^{2}{c_2}\,{c_1}^{4}
+1728\,{c_4}\,{c_3}^{2}{c_2}^{4}{c_1}^{2}
+1728\,{c_4}^{2}{c_3}^{2}{c_2}^{3}\\
&-10800\,{c_4}^{3}{c_3}^{2}{c_2}
-10800\,{c_4}^{4}{c_2}\,{c_1}^{2}
+8640\,{c_4}^{3}{c_2}\,{c_1}^{6}\\
&-864\,{c_4}^{3}{c_3}\,{c_1}^{5}
-12960\,{c_4}^{3}{c_3}^{2}{c_1}^{2}
+1728\,{c_4}^{3}{c_2}^{3}{c_1}^{2}\\
&-1728\,{c_4}^{2}{c_3}^{2}{c_1}^{6}
+7776\,{c_4}^{3}{c_2}^{2}{c_1}^{4}
-11664\,{c_4}\,{c_1}\,{c_3}^{5}\\
&-1728\,{c_4}\,{c_3}^{4}{c_1}^{4}
+32400\,{c_4}^{3}{c_3}\,{c_2}\,{c_1}^{3}
+10800\,{c_4}^{3}{c_3}\,{c_1}\,{c_2}^{2}\\
&-24624\,{c_4}^{4}{c_1}^{4}
-1728\,{c_4}^{2}{c_3}\,{c_1}\,{c_2}^{4}
-12096\,{c_4}^{2}{c_3}\,{c_2}^{2}{c_1}^{5}\\
&-3888\,{c_4}^{2}{c_3}^{3}{c_1}\,{c_2}
-9504\,{c_4}^{2}{c_3}\,{c_2}^{3}{c_1}^{3}
+3456\,{c_4}\,{c_3}^{2}{c_2}^{3}{c_1}^{4}\\
&+19440\,{c_4}\,{c_3}^{4}{c_2}\,{c_1}^{2}
+1728\,{c_4}\,{c_3}^{3}{c_2}\,{c_1}^{5}
-11232\,{c_4}\,{c_3}^{3}{c_2}^{2}{c_1}^{3}
\endaligned{}
$$

$>$ number\_X2:=integral(Gc,cl\_limiting\_P3\_in\_X2):

$>$ latex($''$);

$$
number\_X2 \ = \ 0
$$

\# Construct the projective bundle $\Bbb P({\Sym}^2U^*)$ with
the projection map $\pi_2 = P2$ and the

\# tautological bundle $T_2 = z2$.

$>$ Proj(P2,dual(s2),z2):

\# Use formula (3.3) to compute the class and the number
of limiting $\Bbb P^3$'s in a generic

\# hyperplane.

$>$ cl\_limiting\_P3\_in\_X1:=

$>$ expand(chern(4,Qc)$*$
lowerstar(P2,chern(16,s3-Qc{\&}$*$o(-z2))$*$chern(9,s2-o(-z2)))):

$>$ latex($''$);

$$
\aligned
cl\_limiting\_P3\_in\_X1 \, = \,
&11232\,{c_4}^{2}{c_3}^{3}{c_1}^{3}
-864\,{c_4}\,{c_3}^{3}{c_1}\,{c_2}^{3}
+5832\,{c_4}^{2}{c_3}^{4}\\
&+1296\,{c_4}^{3}{c_2}^{4}
-1728\,{c_4}^{2}{c_3}^{2}{c_2}\,{c_1}^{4}
+864\,{c_4}\,{c_3}^{2}{c_2}^{4}{c_1}^{2}\\
&+50625\,{c_4}^{5}
+864\,{c_4}^{2}{c_3}^{2}{c_2}^{3}
-3780\,{c_4}^{3}{c_3}^{2}{c_2}\\
&+46845\,{c_4}^{4}{c_2}\,{c_1}^{2}
+8856\,{c_4}^{3}{c_2}\,{c_1}^{6}
+34425\,{c_4}^{4}{c_1}\,{c_3}\\
&+1296\,{c_4}^{2}{c_2}^{5}{c_1}^{2}
-16632\,{c_4}^{3}{c_3}\,{c_1}^{5}
-68202\,{c_4}^{3}{c_3}^{2}{c_1}^{2}\\
&-15336\,{c_4}^{3}{c_2}^{3}{c_1}^{2}
+1728\,{c_4}^{2}{c_3}^{2}{c_1}^{6}
+2592\,{c_4}^{2}{c_2}^{4}{c_1}^{4}\\
&-22356\,{c_4}^{3}{c_2}^{2}{c_1}^{4}
-5832\,{c_4}\,{c_1}\,{c_3}^{5}
-864\,{c_4}\,{c_3}^{4}{c_1}^{4}\\
&+55566\,{c_4}^{3}{c_3}\,{c_2}\,{c_1}^{3}
-6912\,{c_4}^{3}{c_3}\,{c_1}\,{c_2}^{2}
-16200\,{c_4}^{4}{c_2}^{2}\\
&+42120\,{c_4}^{4}{c_1}^{4}
+2916\,{c_4}^{2}{c_3}^{2}{c_2}^{2}{c_1}^{2}
+432\,{c_4}^{2}{c_3}\,{c_1}\,{c_2}^{4}\\
&+432\,{c_4}^{2}{c_3}\,{c_2}^{2}{c_1}^{5}
+6804\,{c_4}^{2}{c_3}^{3}{c_1}\,{c_2}
-5400\,{c_4}^{2}{c_3}\,{c_2}^{3}{c_1}^{3}\\
&+1728\,{c_4}\,{c_3}^{2}{c_2}^{3}{c_1}^{4}
+9720\,{c_4}\,{c_3}^{4}{c_2}\,{c_1}^{2}
+864\,{c_4}\,{c_3}^{3}{c_2}\,{c_1}^{5}\\
&-5616\,{c_4}\,{c_3}^{3}{c_2}^{2}{c_1}^{3}
\endaligned{}
$$

$>$ number\_X1:=integral(Gc,cl\_limiting\_P3\_in\_X1):

$>$ latex($''$);

$$
number\_X1 \ = \ 321489
$$

$>$ quit;

{}

\enddocument{}

\bye